\documentclass[12pt]{article}
\usepackage{amssymb}
\usepackage{amsthm}
\usepackage{amsmath}
\usepackage{graphicx}
\usepackage{enumerate}
\newtheorem{theorem}{Theorem}[section]
\newtheorem{definition}[theorem]{Definition}
\newtheorem{lemma}[theorem]{Lemma}
\newtheorem{proposition}[theorem]{Proposition}
\newtheorem{remark}[theorem]{Remark}
\newtheorem{example}[theorem]{Example}
\newtheorem{corollary}[theorem]{Corollary}
\title{Boolean properties and Bell-like inequalities of numerical events}
\author{D.~Dorninger, H.~L\"anger and M.~J.~M\c aczy\'nski}
\date{}
\begin{document}
\footnotetext[1]{Support of the research of the second author by \"OAD, project CZ~02/2019, is gratefully acknowledged.}
\maketitle
\begin{abstract}
Let $S$ be a set of states of a physical system and $p(s)$ be the probability of the occurrence of an event when the system is in state $s\in S$. A function $p:S\rightarrow[0,1]$ is called a numerical event or alternatively an $S$-probability. If a set $P:=\{p(s)\mid s\in S\}$ is ordered by the order of real functions such that certain plausible requirements are fulfilled, $P$ becomes an orthomodular poset in which properties can be described by the addition and comparison of functions. $P$ is then called an algebra of $S$-probabilities or algebra of numerical events. We first answer the question under which circumstances it is possible to consider sets of empirically found numerical events as members of an algebra of $S$-probabilities. Then we discuss the problem to decide whether a given small set $P_n$ of $S$-probabilities can be embedded into a Boolean subalgebra of an algebra $P$ of $S$-probabilities, in which case we will call $P_n$ Boolean embeddable. If $P_n$ is not Boolean embeddable, then the physical system at hand will most likely be non-classical. In the case of a concrete logic $P$, that is a quantum logic which can be represented by sets, we derive criteria for $P_n\subseteq P$ to be Boolean embeddable which can be checked by very simple procedures, for arbitrary $S$-probabilities we provide sets of Bell-like inequalities which characterize the Boolean embeddability of $P_n$. Finally we will show how these Bell-like inequalities fit into a general framework of Bell inequalities by providing a method for generating  Bell inequalities for $S$-probabilities from elementary Bell valuations.
\end{abstract}
 


\section{Introduction}

Let $S$ be a set of states of a physical system and $p(s)$ the probability of an occurrence of an event when the system is in state $s\in S$. The function $p$ from $S$ to $[0,1]$ is called a {\em numerical event}, {\em multidimensional probability} or, more precisely, an {\em $S$-probability} (cf.\ \cite{BM},\cite{BM93}). We note that $p(s)$ can also be considered as a special case of Mackey's probability function $p(A,s,E)$, with $A$ a fixed observable, $s$ a variable state, and $E$ a fixed Borel set (cf.\ \cite M).

\begin{definition}
Let us denote the constant functions $p(s)=0$ and $p(s)=1$ for all $s\in S$ by $0$ and $1$ and write $p\perp q$ if the functions $p$ and $q$ are orthogonal, i.e., $p\leq q':=1-q$. We further agree that the symbols $+$ and $-$ indicate the addition and subtraction in $\mathbb R$. Then a set $P$ of $S$-probabilities is called a {\em space of numerical events} {\rm(}cf.\ {\rm\cite{BM93})} or an {\em algebra of $S$-probabilities} {\rm(}or an {\em algebra of numerical events}{\rm)} {\rm(}first mentioned in {\rm\cite{BDM})}, if it satisfies the following axioms:
\begin{enumerate}
\item [{\rm(A1)}] $0\in P$,
\item [{\rm(A2)}] if $p\in P$ then its {\em complement} {\rm(}considered as a numerical event{\rm)} $p'=1-p\in P$,
\item [{\rm(A3)}] if $p,q\in P$ and $p\perp q$ then $p+q\in P$,
\item [{\rm(A4)}] if $p,q,r\in P$ and $p\perp q\perp r\perp p$ then $p+q+r\in P$.
\end{enumerate}
If only {\rm(A1)} -- {\rm(A3)} are fulfilled, then $P$ is called a {\em generalized field of events}, abbreviated by {\rm GFE} {\rm(}cf.\ {\rm\cite D)}.
\end{definition}

Obviously, under the assumption of (A1), (A3) is a special case of (A4). Moreover, we point out that under the assumption of (A2), Axiom (A3) is equivalent to the following axiom:
\begin{enumerate}
\item[(A5)] If $p,q\in P$ and $p\leq q$ then $q-p\in P$.
\end{enumerate}

If all elements of a GFE $Q$ can only assume the values $0$ and $1$, then $Q$ is an algebra of $S$-probabilities, as one can immediately see, and, moreover (cf.\ \cite{DL09}), $Q$ is a {\em concrete logic}, i.e.\ a logic that has a set representation (cf.\ \cite P). 

Every algebra of $S$-probabilities can be interpreted as the range of a complete probability measure on an orthomodular poset $L$ conceived as an event system (cf.\ \cite{BM}). If $L$ is a classical event system one can show that the associated algebra of $S$-probabilities is a Boolean algebra, and conversely, algebras of $S$-probabilities that are Boolean algebras are the ranges of classical event systems. So it will be crucial to answer the question under which conditions an algebra of $S$-probabilities will be a Boolean algebra, and, in particular, when is this the case if one only knows a few $S$-probabilities obtained by measurements.

Algebras of $S$-probabilities can also be characterized algebraically: Up to isomorphism, they are exactly the orthomodular posets that admit a full set of states (cf.\ \cite{MT}). 

We point out that many well-known quantum logics can be considered as algebras of $S$-probabilities, like all concrete logics and Hilbert space logics: To see this for Hilbert space logics, let $H$ be a Hilbert space, $P(H)$ the set of orthogonal projectors of $H$, $S$ the set of one-dimensional subspaces of $H$, and for every $s\in S$ let $v_s$ be a fixed vector belonging to $s$. Then, denoting the inner product in $H$ by $\langle.,.\rangle$, the set of functions $\{s\mapsto\langle Av_s,v_s\rangle\mid A\in P(H)\}$ is an algebra of $S$-probabilities (cf.\ \cite{BM93}).

First, we will deal with the following question: Given a set $P_n$ of $S$-probabilities obtained by measurements, can $P_n$ be considered as a subset of an appropriate algebra $P$ of $S$-probabilities, and if so, how can one determine such a $P$? We will provide answers depending on the comparability of the numerical events at hand and also in case that the $S$-probabilities can only assume the values $0$ and $1$.

Next we will consider the problem to decide about the classicality and non-classicality of a physical system by means of a series of $S$-probabilities which have been achieved by measurements. In particular, given a set $P_n$ of $n\leq4$ $S$-probabilities ($|S|$ not necessarily limited to a small number) can these numerical events be embedded into a Boolean subalgebra of an algebra of $S$-probabilities? If this is the case we will call $P_n$ {\em Boolean embeddable}. Answering this question, above all one will expect the following insight: If $P_n$ is not Boolean embeddable then one can be quite sure to deal with quantum mechanical effects, i.e., the physical system will not be classical.

For $S$-probabilities which can only assume the values $0$ and $1$ criteria for classicality will be derived which can be most simply checked by elementary operations. For $S$-probabilities that can assume arbitrary values from $[0,1]$ characteristic sets of Bell-like inequalities will be provided. If one of these inequalities is violated for one $s\in S$ the system will not be classical.

Finally, we will consider Bell inequalities from a more general point of view by specifying how to generate Bell inequalities for $S$-probabilities from elementary Bell valuations. We will show that the Bell-like inequalities we set up to characterize the classicality of a physical system do fit into this framework.

\section{Suitability of data}

To be able to utilize the theory of algebras of $S$-probabilities when concrete measurements are available, we first ask whether it is possible in principle to embed empirically found numerical events into an algebra of $S$-probabilities. Of course this will be no problem, if one has enough knowledge about the algebra of $S$-probabilities that governs an experiment, like with the determination of probabilities of the transmission of a photon through a polarizer, that will result in a set of $\{s_1,s_2\}$-probabilities
\[
P=\{(0,0),(1,1)\}\cup\{(\delta,1-\delta)\mid\delta\in[0,1]\setminus\{1/2\}\},
\]
or else, if one carries out measurements within the confines of a known Hilbert space logic pertinent to a physical system. However, if one has no or not sufficient prior knowledge about the underlying algebra of $S$-probabilities then one has to check the data for their suitability as elements of a possibly unknown algebra of $S$-probabilities.

We first observe that the elements $p\neq0,1$ of an algebra of $S$-probabilities have the property (cf.\ \cite{DDL10b}) that neither $p\leq p'$ nor $p\geq p'$. If this is the case we will call $p$ {\em proper}. Obviously, an $S$-probability $\neq0,1$ which can only achieve the values $0$ and $1$ is proper. 
 
Now let 
\[
P_n=\{p_1,\ldots,p_n\}
\]
be a {\em set of $n$ pairwise different proper $S$-probabilities} and
\[
\overline{P_n}:=P_n\cup\{p_1',\ldots,p_n'\},
\]
both partially ordered by the order $\leq$ of functions. We will say that $P_n$ is {\em embeddable} into an algebra $P$ of $S$-probabilities if $P_n$ is order-isomorphic to a subset of $P$.

\begin{proposition}\label{Prop1}
If every $p\in P_n$ can only assume the values $0$ and $1$ then $P_n$ is always embeddable into an algebra of $S$-probabilities.
\end{proposition}

\begin{proof}
One can consider the elements $0$ and $1$ of a GFE as nullary operations, $'$ as a unary operation and $+$ as a partial operation. Then the structure generated by the elements of $P_n$ by means of the operations $0$, $1$, $'$ and $+$ is a GFE $[P_n]_{01}$ (cf.\ \cite D) which, because of its two-valued functions, is also an algebra of $S$-probabilities. Moreover, the order of the elements of $P_n$ will not be altered within $[P_n]_{01}$.
\end{proof}

As it is very common within the theory of lattices we will denote the poset consisting of $0$ and $1$ and an ``antichain'' (pairwise different and incomparable elements) of elements $a_1,\ldots,a_n,a_1',\ldots,a_n'$ between $0$ and $1$ by ${\rm MO}_n$.

\begin{proposition}\label{Prop2}
Assume $n=2$.
\begin{enumerate}[{\rm(a)}]
\item If $\overline{P_2}$ is a four-element antichain then there do exist algebras of $S$-probabilities into which $P_2$ can be embedded. The smallest such algebra of $S$-probabilities is ${\rm MO}_2$ and the smallest Boolean algebra has 16 elements.
\item If $\overline{P_2}$ is not a four-element antichain, say $p_1<p_2$, and if $p_2-p_1$ is proper, then there do exist algebras of $S$-probabilities into which $P_2$ can be embedded. An eight-element Boolean algebra is the smallest algebra of $S$-probabilities to contain $P_2$.
\end{enumerate}
\end{proposition}

\begin{proof}
\
\begin{enumerate}[(a)]
\item Obviously the orthomodular lattice MO$_2$ contains $P_2$, and, as a straightforward calculation shows, the smallest Boolean algebra to comprise $P_2$ as a suborder will be an algebra of 16 elements.
\item If we assume without loss of generality that $p_1<p_2$ and we replenish $p_1,p_2,p_1',p_2'$ by $p_{12}:=p_2-p_1,p_{12}',0,1$, respectively, then we obtain an eight-element Boolean algebra which contains $P_2$, and no smaller algebra of $S$-probabilities will exist to comprise $P_2$.
\end{enumerate}
\end{proof}

The first statement of Proposition~\ref{Prop2} can be immediately generalized in the following way:  

\begin{proposition}\label{Prop3}
$P_n$ can be embedded into an algebra of $S$-probabilities if all elements of $\overline{P_n}$ are pairwise incomparable, and ${\rm MO}_n$ is the smallest algebra of $S$-probabilities to encompass $P_n$.
\end{proposition}
 
\begin{proposition}\label{Prop4}
Assume $|S|=2$. Then $P_n$ is embeddable into an algebra of $S$-prob\-a\-bil\-i\-ties if and only if the elements of $\overline{P_n}$ are pairwise incomparable.
\end{proposition}

\begin{proof}
If $p:=(a_1,a_2)\leq(b_1,b_2)=:q$ and $a_1<1/2$ then $a_2>1/2$, $b_2>1/2$, $b_1<1/2$. Hence $b_1-a_1<1/2$, $b_2-a_2<1/2$ which means that $q-p$ is not proper contrary to the fact (see above) that $q-p$ has to be an element of the supposed algebra of $S$-probabilities. In case $a_1>1/2$ one can conclude along the same lines that $q-p$ is not proper.
\end{proof}

Supposed that one does not have a sufficient knowledge of the algebra of $S$-probabilities to which $P_n$ belongs, then Propositions~\ref{Prop1} -- \ref{Prop4} illustrate that it is only possible to utilize the theory of $S$-probabilities for an arbitrary set of $S$-probabilities gained by measurements, if all numerical events have only the outcome $0$ and $1$, in other words, if one deals with a concrete logic. In all other cases only special sets of data will be of use.

\section{Concrete logics of numerical events}

To begin with let $P$ be an arbitrary algebra of $S$-probabilities and $p,q\in P$. If the {\em infimum} of $p$ and $q$ exists in $P$ we will denote it by $p\wedge q$, if their {\em supremum} exists in $P$ we will denote it by $p\vee q$. Moreover, we point out that if $p\perp q$ then $p+q=p\vee q$, and if $p\leq q$ then $q-p\in P$ (see (A5)) and $q-p=q\wedge p'$ (cf.\ \cite{DDL10b}).

Two elements $f$ and $g$ are said to {\em commute} -- in symbols $f\mathrel{{\rm C}}g$ -- if there exists some Boolean subalgebra $B$ of $P$ that contains $f$ and $g$. Here Boolean subalgebra means that the operations $0$, $1$, complementation, forming of infima and suprema in $B$ coincide with $0$, $1$, $'$, $\wedge$, $\vee$ in $P$, respectively.

As shown in \cite{DL14b}, $f\mathrel{{\rm C}}g$ holds in $P$ if and only if there exists some $a\in P$ with $a\leq f\leq a+g\leq1$, in which case we write $f\mathrel{{\rm C}(a)}g$. If $f$ and $g$ commute, so do $f$ and $g'$, $f'$ and $g$ as well as $f'$ and $g'$. Two elements that are comparable always commute.

Now assume that $P$ is a {\em concrete logic} which means that the range of the elements of $P$ is $\{0,1\}$.

\begin{definition}
For $f,g\in P$ let $f\barwedge g$ be the $S$-probability $h$ with $h(s)=\min\{f(s),g(s)\}$ for $s\in S$.
\end{definition}

\begin{lemma}\label{lem1}
Let $P$ be a concrete logic and $f,g,a,b\in P$. Then $f\mathrel{{\rm C}(a)}g$ if and only if $a=f\wedge g'=f\barwedge g'$. Hence, we have $f\mathrel{{\rm C}(b)}g'$ if and only if $b=f\wedge g=f\barwedge g$.
\end{lemma}

\begin{proof}
According to the definition of C$(a)$, $f\mathrel{{\rm C}(a)}g$ is equivalent to $a\leq f\leq a+g\leq1$. Now consider $s\in S$. Then one can easily see that in all four cases for $(f(s),g(s))$ the value of $a(s)$ is uniquely determined and coincides with $\min\{f(s),g'(s)\}$.
\[
\begin{array}{c|c|c}
f(s) & g(s) & a(s) \\
\hline
 0   &  0   &  0 \\
 0   &  1   &  0 \\
 1   &  0   &  1 \\
 1   &  1   &  0
\end{array}
\]
This verifies $a=f\wedge g'=f\barwedge g'$.
\end{proof}

For the sake of a simpler lettering we agree to the following understanding within proofs of theorems of this section, justified by Lemma \ref{lem1}: For elements $a=f\barwedge g'$ that are used in connection with a relation $f\mathrel{{\rm C}(a)}g$ we will write $f\wedge g'$ instead of $f\barwedge g'$, which means that if an infimum will be considered in combination with a relation $C(a)$ we will omit the bar in $\barwedge$. Since $f\mathrel{{\rm C}}g$ is equivalent to $f\mathrel{{\rm C}}g'$, $f'\mathrel{{\rm C}}g$ and $f'\mathrel{{\rm C}}g'$ we will also use this convention with $f\wedge g$, $f\wedge g'$, $f'\wedge g$ and $f'\wedge g'$, and by de Morgan's laws we can also agree to proceed the same way with suprema $f\,\overline{\vee}\,g:=(f'\barwedge g')'$. Of course, in general $f\wedge g$ does not coincide with $f\barwedge g$, and the same is true for suprema.

\begin{remark}\label{rem1}
For $f,g\in P$ the following assertions are equivalent:
\begin{itemize}
\item There exists some $a\in P$ with $f\mathrel{{\rm C}(a)}g$,
\item there exists some $b\in P$ with $f\mathrel{{\rm C}(b)}g'$.
\end{itemize}
\end{remark}

If a set $P_n=\{p_1,\ldots,p_n\}$ of $n\geq 2$ elements is contained in a Boolean subalgebra of an \textit{arbitrary} algebra of $S$-probabilities $P$ we have called it Boolean embeddable. In accordance with \cite{BM93} and \cite{DL14b} we will now say for the sake of brevity that $P_n$ is {\em Boolean} (though the notion of Boolean is also used for posets with a different meaning). If $P_n$ is Boolean, measurements resulting in $P_n$ will indicate that one deals with a classical physical system.

In the sequel we will make frequent use of the following three theorems from \cite{DL14b} holding for arbitrary algebras of $S$-probabilities. 

\begin{theorem}\label{th1}
$P_2=\{p_1,p_2\}$ is Boolean if and only if there exists some
\[
a_{12}\in P
\]
satisfying
\begin{align*}
p_1 & \mathrel{{\rm C}(a_{12})}p_2.
\end{align*} 
\end{theorem}

\begin{theorem}\label{th2}
$P_3=\{p_1,p_2,p_3\}$ is Boolean if and only if there exist
\[
a_{12},a_{13},a_{23},a_{1213}\in P
\]
satisfying
\begin{align*}
         p_i & \mathrel{{\rm C}(a_{ij})}p_j\text{ for }i,j\in\{1,2,3\}\text{ with }i<j, \\
(p_1-a_{12}) & \mathrel{{\rm C}(a_{1213})}(p_i-a_{i3})\text{ for }i\in\{1,2\}.
\end{align*}
\end{theorem}

\begin{theorem}\label{th3}
$P_4=\{p_1,\ldots,p_4\}$ is Boolean if and only if there exist
\[
a_{12},a_{13},a_{14},a_{23},a_{24},a_{34},a_{1213},a_{1214},a_{1314},a_{2324},a_{1234},a_{1324},a_{1423},a_{123124},a_{124134},a_{134234}\in P
\]
satisfying
\begin{align*}
                    p_i & \mathrel{{\rm C}(a_{ij})}p_j\text{ for }i,j\in\{1,2,3,4\}\text{ with }i<j, \\
           (p_i-a_{ij}) & \mathrel{{\rm C}(a_{ijik})}(p_i-a_{ik})\text{ for }i,j,k\in\{1,2,3,4\}\text{ with }i<j<k, \\
           (p_i-a_{ij}) & \mathrel{{\rm C}(a_{ijik})}(p_j-a_{jk})\text{ for }i,j,k\in\{1,2,3,4\}\text{ with }i<j<k, \\
           (p_1-a_{12}) & \mathrel{{\rm C}(a_{1234})}(p_3-a_{34}), \\
           (p_1-a_{13}) & \mathrel{{\rm C}(a_{1324})}(p_2-a_{24}), \\
           (p_1-a_{14}) & \mathrel{{\rm C}(a_{1423})}(p_2-a_{23}), \\
((p_1-a_{12})-a_{1213}) & \mathrel{{\rm C}(a_{123124})}((p_i-a_{ij})-a_{iji4})\text{ for }i,j\in\{1,2,3\}\text{ with }i<j, \\
((p_1-a_{12})-a_{1214}) & \mathrel{{\rm C}(a_{124134})}((p_i-a_{i3})-a_{i3i4})\text{ for }i\in\{1,2\}, \\
((p_1-a_{13})-a_{1314}) & \mathrel{{\rm C}(a_{134234})}((p_2-a_{23})-a_{2324}).
\end{align*}
\end{theorem}

Now we assume again that all elements will belong to a concrete logic $P$ and $P_n=\{p_1,\ldots,p_n\}$ is a set of $n\geq 2$ pairwise different elements of $P$.
 
\begin{theorem}
$P_2=\{p_1,p_2\}$ is Boolean if and only if
\[
p_1\barwedge p_2\in P.
\]
\end{theorem}

\begin{proof}
We use Theorem~\ref{th1}. The following three claims are equivalent:
\begin{align*}
& \text{There exists some }a_{12}\in P\text{ with }p_1\mathrel{{\rm C}(a_{12})}p_2, \\
& \text{there exists some }a\in P\text{ with }p_1\mathrel{{\rm C}(a)}p_2', \\
& p_1\wedge p_2\in P.
\end{align*}
\end{proof}

\begin{theorem}\label{theo3}
$P_3=\{p_1,p_2,p_3\}$ is Boolean if and only if
\[
p_i\barwedge p_j\barwedge p_k\in P
\]
for all $i,j,k\in\{1,2,3\}$.
\end{theorem}

\begin{proof}
Now we use Theorem~\ref{th2}. For $i,j\in\{1,2,3\}$ with $i<j$ the following assertions are equivalent:
\begin{align*}
& \text{There exists some }a_{ij}\in P\text{ with }p_i\mathrel{{\rm C}(a_{ij})}p_j, \\
& \text{there exists some }b_{ij}\in P\text{ with }p_i\mathrel{{\rm C}(b_{ij})}p_j', \\
& p_i\wedge p_j\in P.
\end{align*}
And for $i,j\in\{1,2,3\}$ with $i<j$ we have
\begin{align*}
    a_{ij} & =p_i\wedge p_j', \\
p_i-a_{ij} & =p_i\wedge(p_i\wedge p_j')'=p_i\wedge(p_i'\vee p_j)=p_i\wedge p_j.
\end{align*}
Therefore,
\[
(p_1-a_{12})\wedge(p_i-a_{i3})'=(p_1\wedge p_2)\wedge(p_i\wedge p_3)'=(p_1\wedge p_2)\wedge(p_i'\vee p_3')=p_1\wedge p_2\wedge p_3'
\]
for $i\in\{1,2\}$, from which can we infer that the following four claims are equivalent:
\begin{align*}
& \text{There exists some }a_{1213}\in P\text{ with }(p_1-a_{12})\mathrel{{\rm C}(a_{1213})}(p_i-a_{i3})\text{ for }i\in\{1,2\}, \\
& \text{there exists some }a\in P\text{ with }(p_1-a_{12})\mathrel{{\rm C}(a)}(p_1-a_{13}), \\
& \text{there exists some }b\in P\text{ with }(p_1-a_{12})\mathrel{{\rm C}(b)}(p_1-a_{13})', \\
& p_1\wedge p_2\wedge p_3\in P.
\end{align*}
\end{proof}

\begin{theorem}\label{theo4}
$P_4=\{p_1,\ldots,p_4\}$ is Boolean if and only if
\[
p_i\barwedge p_j\barwedge p_k\barwedge p_m\in P
\]
for all $i,j,k,m\in\{1,\ldots,4\}$.
\end{theorem}

\begin{proof}
Here we apply Theorem~\ref{th3}. For $i,j\in\{1,\ldots,4\}$ with $i<j$ 
\[
\text{there exists some }a_{ij}\in P\text{ with }p_i\mathrel{{\rm C}(a_{ij})}p_j
\]
if and only if
\[
\text{there exists some }b_{ij}\in P\text{ with }p_i\mathrel{{\rm C}(b_{ij})}p_j'
\]
which is equivalent to
\[
p_i\wedge p_j\in P.
\]
For $i,j\in\{1,\ldots,4\}$ with $i<j$ 
\begin{align*}
    a_{ij} & =p_i\wedge p_j', \\
p_i-a_{ij} & =p_i\wedge(p_i\wedge p_j')'=p_i\wedge(p_i'\vee p_j)=p_i\wedge p_j.
\end{align*}
Further, for $i,j,k\in\{1,\ldots,4\}$ with $i<j<k$ 
\begin{align*}
& (p_i-a_{ij})\wedge(p_i-a_{ik})'=(p_i\wedge p_j)\wedge(p_i\wedge p_k)'=(p_i\wedge p_j)\wedge(p_i'\vee p_k')=p_i\wedge p_j\wedge p_k', \\
& (p_i-a_{ij})\wedge(p_j-a_{jk})'=(p_i\wedge p_j)\wedge(p_j\wedge p_k)'=(p_i\wedge p_j)\wedge(p_j'\vee p_k')=p_i\wedge p_j\wedge p_k'.
\end{align*}
Therefore the following assertions are equivalent for $i,j,k\in\{1,\ldots,4\}$ with $i<j<k$:
\begin{align*}
& \text{There exists some }a_{ijik}\in P\text{ with }(p_i-a_{ij})\mathrel{{\rm C}(a_{ijik})}(p_i-a_{ik})\text{ and} \\
& \quad\quad(p_i-a_{jk})\mathrel{{\rm C}(a_{ijik})}(p_j-a_{jk}), \\
& \text{there exists some }a\in P\text{ with }(p_i-a_{ij})\mathrel{{\rm C}(a)}(p_i-a_{ik}), \\
& \text{there exists some }b\in P\text{ with }(p_i-a_{ij})\mathrel{{\rm C}(b)}(p_i-a_{ik})', \\
& p_i\wedge p_j\wedge p_k\in P.
\end{align*}
Next we take into account that the following three statements are equivalent:
\begin{align*}
& \text{There exists some }a_{1234}\in P\text{ with }(p_1-a_{12})\mathrel{{\rm C}(a_{1234})}(p_3-a_{34}), \\
& \text{there exists some }c\in P\text{ with }(p_1-a_{12})\mathrel{{\rm C}(c)}(p_3-a_{34})', \\
& p_1\wedge\cdots\wedge p_4\in P.
\end{align*}
Also equivalent are:
\begin{align*}
& \text{There exists some }a_{1324}\in P\text{ with }(p_1-a_{13})\mathrel{{\rm C}(a_{1324})}(p_2-a_{24}), \\
& \text{there exists some }d\in P\text{ with }(p_1-a_{13})\mathrel{{\rm C}(d)}(p_2-a_{24})', \\
& p_1\wedge\cdots\wedge p_4\in P.
\end{align*}
The same is true for the statements
\begin{align*}
& \text{There exists some }a_{1423}\in P\text{ with }(p_1-a_{14})\mathrel{{\rm C}(a_{1423})}(p_2-a_{23}), \\
& \text{there exists some }e\in P\text{ with }(p_1-a_{14})\mathrel{{\rm C}(e)}(p_2-a_{23})', \\
& p_1\wedge\cdots\wedge p_4\in P.
\end{align*}
For $i,j,k\in\{1,\ldots,4\}$ with $i<j<k$ we have
\begin{align*}
             a_{ijik} & =p_i\wedge p_j\wedge p_k',\text{ therefore} \\
(p_i-a_{ij})-a_{ijik} & =(p_i\wedge p_j)\wedge(p_i\wedge p_j\wedge p_k')'=(p_i\wedge p_j)\wedge(p_i'\vee p_j'\vee p_k)=p_i\wedge p_j\wedge p_k.
\end{align*}
As for $i,j\in\{1,2,3\}$ with $i<j$
\[
\text{there exists some }a_{123124}\in P\text{ with }((p_1-a_{12})-a_{1213})\mathrel{{\rm C}(a_{123124})}((p_i-a_{ij})-a_{iji4})
\]
if and only if
\[
\text{there exists some }f\in P\text{ with }((p_1-a_{12})-a_{1213})\mathrel{{\rm C}(f)}((p_i-a_{ij})-a_{iji4})'
\]
which is equivalent to
\[
p_1\wedge\cdots\wedge p_4\in P.
\]
Moreover, for $i\in\{1,2\}$ the following three assertions are equivalent:
\begin{align*}
& \text{There exists some }a_{124134}\in P\text{ with }((p_1-a_{12})-a_{1214})\mathrel{{\rm C}(a_{124134})}((p_i-a_{i3})-a_{i3i4}), \\
& \text{there exists some }g\in P\text{ with }((p_1-a_{12})-a_{1214})\mathrel{{\rm C}(g)}((p_i-a_{i3})-a_{i3i4})', \\
& p_1\wedge\cdots\wedge p_4\in P.
\end{align*}
Finally, we take into account the equivalence of the assertions:
\begin{align*}
& \text{There exists some }a_{134234}\in P\text{ with }((p_1-a_{13})-a_{1314})\mathrel{{\rm C}(a_{134234})}((p_2-a_{23})-a_{2324}), \\
& \text{there exists some }h\in P\text{ with }((p_1-a_{13})-a_{1314})\mathrel{{\rm C}(h)}((p_2-a_{23})-a_{2324})', \\
& p_1\wedge\cdots\wedge p_4\in P.
\end{align*}
\end{proof}

\section{Bell-like inequalities of numerical events}

Now we do not any more restrict the values of $S$-probabilities to be only $0$ and $1$.

\begin{definition}
Let $P$ be an arbitrary algebra of $S$-probabilities and $P_n$ a set of pairwise different proper $S$-probabilities of $P$. We call $p_i,p_j\in P$ correlated in $P$ if the element $p_{ij}:=p_i\wedge p_j$ exists in $P$. Moreover, we call $p_i,p_j,p_k$ correlated, if the element $p_{ijk}:=p_i\wedge p_j\wedge p_k$ exists in $P$, and so forth, for more than three elements.
\end{definition}

Now we turn to the question whether a set $P_n$ which contains correlated elements will be Boolean.

\begin{theorem}\label{P2}
Assume that the elements $p_1,p_2$ of $P_2$ are correlated with $p_{12}=p_1\wedge p_2$. Then $P_2$ is Boolean if and only if
\begin{align*}
p_1+p_2-p_{12} & \leq1.
\end{align*}
\end{theorem}

\begin{proof}
If $P_2$ is Boolean then
\begin{align*}
p_1+p_2-p_{12} & =p_1\vee p_2\leq1.
\end{align*}
If, conversely, the inequality of the theorem holds, then
\begin{align*}
a_{12} & :=p_1-p_{12}
\end{align*}
satisfies the condition of Theorem~\ref{th1}.
\end{proof}

\begin{theorem}\label{P3}
Assume that the elements $p_1,p_2,p_3$ of $P_3$ are correlated with $p_{ij}=p_i\wedge p_j$ for $i,j\in\{1,2,3\}$ with $i<j$, and $p_{123}=p_1\wedge p_2\wedge p_3$. Then $P_3$ is Boolean if and only if
\begin{align*}
       p_i+p_j-p_{ij} & \leq1\text{ for }i,j\in\{1,2,3\}\text{ with }i<j, \\
p_{12}+p_{i3}-p_{123} & \leq1\text{ for }i\in\{1,2\}.
\end{align*}
\end{theorem}

\begin{proof}
If $P_3$ is Boolean then
\begin{align*}
       p_i+p_j-p_{ij} & =p_i\vee p_j\leq1\text{ for }i,j\in\{1,2,3\}\text{ with }i<j, \\
p_{12}+p_{i3}-p_{123} & =p_{12}\vee p_{i3}\leq1\text{ for }i\in\{1,2\}.
\end{align*}
If, conversely, the inequalities of the theorem hold then
\begin{align*}
  a_{ij} & :=p_i-p_{ij}\text{ for }i,j\in\{1,2,3\}\text{ with }i<j, \\
a_{1213} & :=p_{12}-p_{123}
\end{align*}
meet the requirements of Theorem~\ref{th2}.
\end{proof}

\begin{theorem}\label{P4}
Assume that the elements $p_1,\ldots,p_4$ of $P_4$ are correlated with $p_{ij}=p_i\wedge p_j$ for $i,j\in\{1,\ldots,4\}$ with $i<j$, $p_{ijk}=p_i\wedge p_j\wedge p_k$ for $i,j,k\in\{1,\ldots,4\}$ with $i<j<k$ and $p_{1234}=p_1\wedge\cdots\wedge p_4$. Then $P_4$ is Boolean if and only if
\begin{align*}
          p_i+p_j-p_{ij} & \leq1\text{ for }i,j\in\{1,\ldots,4\}\text{ with }i<j, \\
   p_{ij}+p_{ik}-p_{ijk} & \leq1\text{ for }i,j,k\in\{1,\ldots,4\}\text{ with }i<j<k, \\
   p_{ij}+p_{jk}-p_{ijk} & \leq1\text{ for }i,j,k\in\{1,\ldots,4\}\text{ with }i<j<k, \\
  p_{12}+p_{34}-p_{1234} & \leq1, \\
  p_{13}+p_{24}-p_{1234} & \leq1, \\
  p_{14}+p_{23}-p_{1234} & \leq1, \\
p_{123}+p_{ij4}-p_{1234} & \leq1\text{ for }i,j\in\{1,2,3\}\text{ with }i<j, \\
p_{124}+p_{i34}-p_{1234} & \leq1\text{ for }i\in\{1,2\}, \\
p_{134}+p_{234}-p_{1234} & \leq1.
\end{align*}
\end{theorem}

\begin{proof}
If $P_4$ is Boolean then
\begin{align*}
          p_i+p_j-p_{ij} & =p_i\vee p_j\leq1\text{ for }i,j\in\{1,\ldots,4\}\text{ with }i<j, \\
   p_{ij}+p_{ik}-p_{ijk} & =p_{ij}\vee p_{ik}\leq1\text{ for }i,j,k\in\{1,\ldots,4\}\text{ with }i<j<k, \\
   p_{ij}+p_{jk}-p_{ijk} & =p_{ij}\vee p_{jk}\leq1\text{ for }i,j,k\in\{1,\ldots,4\}\text{ with }i<j<k, \\
  p_{12}+p_{34}-p_{1234} & =p_{12}\vee p_{34}\leq1, \\
  p_{13}+p_{24}-p_{1234} & =p_{13}\vee p_{24}\leq1, \\
  p_{14}+p_{23}-p_{1234} & =p_{14}\vee p_{23}\leq1, \\
p_{123}+p_{ij4}-p_{1234} & =p_{123}\vee p_{ij4}\leq1\text{ for }i,j\in\{1,2,3\}\text{ with }i<j, \\
p_{124}+p_{i34}-p_{1234} & =p_{124}\vee p_{i34}\leq1\text{ for }i\in\{1,2\}, \\
p_{134}+p_{234}-p_{1234} & =p_{134}\vee p_{234}\leq1.
\end{align*}
If, conversely, the inequalities of the theorem hold then
\begin{align*}
    a_{ij} & :=p_i-p_{ij}\text{ for }i,j\in\{1,\ldots,4\}\text{ with }i<j, \\
  a_{ijik} & :=p_{ij}-p_{ijk}\text{ for }i,j,k\in\{1,\ldots,4\}\text{ with }i<j<k, \\
  a_{1234} & :=p_{12}-p_{1234}, \\
  a_{1324} & :=p_{13}-p_{1234}, \\
  a_{1423} & :=p_{14}-p_{1234}, \\
a_{123124} & :=p_{123}-p_{1234}, \\
a_{124134} & :=p_{124}-p_{1234}, \\
a_{134234} & :=p_{134}-p_{1234}
\end{align*}
fulfil the conditions of Theorem~\ref{th3}.
\end{proof}

\section{Bell-type inequalities generated by Bell valuations}

The original Bell inequalities and Clauser-Horne-Shimony-Holt inequalities (CHSH-in\-e\-qua\-li\-ties) have been generalized a number of times. Here we refer to the generalizations and, in particular, to the notion of Bell-type inequalities studied in \cite B and \cite{BM} -- \cite{BM94}. We will show how one can generate Bell-type inequalities on the elements from $P$ by using elementary Bell valuations and then classify the above introduced Bell-like inequalities within this context.

Let $P$ be a set of $S$-probabilities which is a Boolean algebra and let $n$ be a positive integer. Further put $N:=\{1,\ldots,n\}$.

\begin{definition}
A function $f$ from $2^N\setminus\{\emptyset\}$ to $\mathbb R$ satisfying
\[
0\leq\sum_{\emptyset\neq J\subseteq I}f(J)\leq1\text{ for all }I\in2^N\setminus\{\emptyset\}
\]
is called a {\em Bell valuation} on $2^N$. The set of all Bell valuations on $2^N$ will be denoted by $B_n$.
\end{definition}

\begin{theorem}\label{th4}
Let $f$ be a Bell valuation on $2^N$. Then the following inequality holds:
\begin{equation}\label{equ2}
0\leq\sum_{\emptyset\neq I\subseteq N}f(I)p_I\leq1
\end{equation}
where $p_1,\ldots,p_n\in P$ and $p_I=\bigwedge\limits_{i\in I}p_i$ for all $I\in2^N\setminus\{\emptyset\}$. Inequality {\rm(\ref{equ2})} is called a \em{Bell-type inequality}.
\end{theorem}

The proof of Theorem~\ref{th4} is given in \cite{BM93} and repeated in \cite B. Theorem~\ref{th4} shows that every Bell valuation generates a Bell-type inequality. 

In the sequel we will explain how to construct Bell valuations from so-called elementary Bell valuations.

\begin{definition}\label{def1}
For every $I\in2^N\setminus\{\emptyset\}$ the function $f_I$ from $2^N\setminus\{\emptyset\}$ to $\mathbb R$ defined by
\[
f_I(J):=\left\{
\begin{array}{ll}
(-1)^{|J\setminus I|} & \text{if }J\supseteq I \\
0                     & \text{otherwise}
\end{array}
\right.
\]
{\rm(}$J\in2^N\setminus\{\emptyset\}${\rm)} is called an {\em elementary Bell valuation} on $2^N$.
\end{definition}

We want to show that every elementary Bell valuation on $2^N$ is a Bell valuation on $2^N$. In fact, we prove more, but first we need some auxiliary results.

\begin{lemma}
Let $I$ be a finite non-empty set. Then half of its subsets have an even number of elements and half an odd number.
\end{lemma}

\begin{proof}
Let $a\in I$ and define a mapping $f$ from $2^I$ to $2^I$ by
\[
f(J):=\left\{
\begin{array}{ll}
J\setminus\{a\} & \text{if }a\in J \\
J\cup\{a\}      & \text{otherwise}
\end{array}
\right.
\]
($J\in2^I$). Then $f$ is an involution (i.e.\ $f\circ f$ is the identical mapping) and hence bijective. Obviously, $f$ maps subsets of $I$ with an even number of elements onto subsets of $I$ with an odd number of elements and vice versa.
\end{proof}

\begin{corollary}
Let $I$ be a finite set. Then
\[
\sum_{J\subseteq I}(-1)^{|J|}=\delta_{I,\emptyset}
\]
where
\[
\delta_{J,K}=\left\{
\begin{array}{rl}
1 & \text{if }J=K \\
0 & \text{otherwise.}
\end{array}
\right.
\]
\end{corollary}

\begin{definition}
For every $h\in\mathbb R^{2^N\setminus\{\emptyset\}}$ let $f_h$ and $g_h$ denote the mappings from $2^N\setminus\{\emptyset\}$ to $\mathbb R$ defined by
\begin{align*}
f_h(I) & :=\sum_{\emptyset\neq J\subseteq N}h(J)f_J(I), \\
g_h(I) & :=\sum_{\emptyset\neq J\subseteq I}h(J)
\end{align*}
for all $I\in2^N\setminus\{\emptyset\}$.
\end{definition}

\begin{lemma}
The mappings $h\mapsto f_h$ and $h\mapsto g_h$ are mutually inverse bijections between $\mathbb R^{2^N\setminus\{\emptyset\}}$ and $\mathbb R^{2^N\setminus\{\emptyset\}}$.
\end{lemma}

\begin{proof}
Obviously, both mappings are well-defined functions from $\mathbb R^{2^N\setminus\{\emptyset\}}$ to $\mathbb R^{2^N\setminus\{\emptyset\}}$. Moreover, if $h\in\mathbb R^{2^N\setminus\{\emptyset\}}$ then
\begin{align*}
g_{f_h}(I) & =\sum_{\emptyset\neq J\subseteq I}\sum_{\emptyset\neq K\subseteq J}h(K)(-1)^{|J\setminus K|}=\sum_{\emptyset\neq K\subseteq I}h(K)\sum_{K\subseteq J\subseteq I}(-1)^{|J\setminus K|}= \\
           & =\sum_{\emptyset\neq K\subseteq I}h(K)\sum_{L\subseteq I\setminus K}(-1)^{|L|}=\sum_{\emptyset\neq K\subseteq I}h(K)\delta_{I\setminus K,\emptyset}=h(I), \\
f_{g_h}(I) & =\sum_{\emptyset\neq J\subseteq I}(\sum_{\emptyset\neq K\subseteq J}h(K))(-1)^{|I\setminus J|}=\sum_{\emptyset\neq K\subseteq I}h(K)\sum_{K\subseteq J\subseteq I}(-1)^{|I\setminus J|}= \\
           & =\sum_{\emptyset\neq K\subseteq I}h(K)\sum_{L\subseteq I\setminus K}(-1)^{|L|}=\sum_{\emptyset\neq K\subseteq I}h(K)\delta_{I\setminus K,\emptyset}=h(I)
\end{align*}
for all $I\in2^N\setminus\{\emptyset\}$.
\end{proof}

\begin{corollary}\label{cor1}
\
\begin{itemize}
\item Let $f\in\mathbb R^{2^N\setminus\{\emptyset\}}$. Then
\begin{itemize}
\item $f(2^N\setminus\{\emptyset\})\subseteq\mathbb Z$ if and only if $g_f(2^N\setminus\{\emptyset\})\subseteq\mathbb Z$,
\item $f\in B_n$ if and only if $g_f(2^N\setminus\{\emptyset\})\subseteq[0,1]$.
\end{itemize}
\item Let $g\in\mathbb R^{2^N\setminus\{\emptyset\}}$. Then
\begin{itemize}
\item $f_g(2^N\setminus\{\emptyset\})\subseteq\mathbb Z$ if and only if $g(2^N\setminus\{\emptyset\})\subseteq\mathbb Z$,
\item $f_g\in B_n$ if and only if $g(2^N\setminus\{\emptyset\})\subseteq[0,1]$.
\end{itemize}
\item We have $f_I=f_{\delta_{JI}}\in B_n$ for all $I\in2^N\setminus\{\emptyset\}$.
\end{itemize}
\end{corollary}

As for applications in quantum mechanics, the original Bell inequalities and the CHSH-inequalities involve only at most four simultaneous events, and only coefficients $1$ (and implicitly $0$) turn up in the inequalities. Therefore linear combinations of elementary Bell valuations with coefficients $0$ and $1$ only are of special importance. It is about Bell-type inequalities of the form $0\leq L\leq1$, where $L$ is a linear combination with coefficients from $\mathbb Z$ of the $S$-probabilities $p_1,\ldots,p_n$ and correlations $p_I=\bigwedge\limits_{i\in I}p_i$, $\emptyset\neq I\subseteq N$. As we have already emphasized a Bell-type inequality is valid if it can be derived in a classical probability space, i.e.\ in a Boolean algebra. By a theorem in \cite{BM93}, a Bell-type inequality
\[
0\leq\sum_{\emptyset\neq I\subseteq N}f(I)p_I\leq1
\]
is valid if and only if $f\in B_n$.

According to Corollary~\ref{cor1}, there exist $2^{2^n-1}-1$ non-zero Bell valuations $f$ on $2^N$ with $f(2^N\setminus\{\emptyset\})\subseteq\mathbb Z$. For $n=2$ this number equals $7$ and for $n=3$ it equals $127$.

\begin{example}
If $n=3$ and
\[
f(I)=\left\{
\begin{array}{rl}
 1 & \text{if }I=\{1\}\text{ or }I=\{2,3\} \\
-1 & \text{if }I=\{1,2\}\text{ or }I=\{1,3\} \\
 0 & \text{otherwise}
\end{array}
\right.
\]
{\rm(}$I\in2^N\setminus\{\emptyset\}${\rm)} then
\[
g_f(I)=\left\{
\begin{array}{rl}
1 & \text{if }I=\{1\}\text{ or }I=\{2,3\} \\
0 & \text{otherwise}
\end{array}
\right.
\]
{\rm(}$I\in2^N\setminus\{\emptyset\}${\rm)} and hence $f\in B_3$ according to Corollary~\ref{cor1}.
\end{example}

\begin{example}
We have
\[
f_{1-\delta_{I,N}}(J)=\sum_{\emptyset\neq K\subseteq J}(-1)^{|J\setminus K|}=\sum_{L\subsetneqq J}(-1)^{|L|}=-(-1)^{|J|}=(-1)^{|J|+1}
\]
for all $J\in2^N\setminus\{\emptyset,N\}$ and
\[
f_{1-\delta_{I,N}}(N)=\sum_{\emptyset\neq J\subseteq N}(1-\delta_{J,N})(-1)^{|N\setminus J|}=\sum_{\emptyset\neq J\subsetneqq N}(-1)^{|N\setminus J|}=\sum_{\emptyset\neq K\subsetneqq N}(-1)^{|K|}=-1-(-1)^n.
\]
\end{example}

Bell inequalities can be used to test classicality or non-classicality of probability systems. If we have a system of $n$ $S$-probabilities $p_1,\ldots,p_n$ and correlations $p_I,\emptyset\neq I\subseteq N,$ then if the system is classical, these probabilities and correlations should satisfy all Bell-type inequalities generated by Bell valuations on $2^N$, in particular the Bell-type inequality generated by the sum of all elementary Bell valuations on $2^N$. This inequality is the following:
\[
0\leq\sum_{\emptyset\neq I\subseteq N}(-1)^{|I|+1}p_I\leq1
\]
since
\[
f_1(I)=\sum_{\emptyset\neq J\subseteq I}(-1)^{|I\setminus J|}=\sum_{K\subsetneqq I}(-1)^{|K|}=-(-1)^{|I|}=(-1)^{|I|+1}
\]
for all $I\in2^N\setminus\{\emptyset\}$. If this inequality is violated, the system is not classical. 

For $n=3$ this inequality has the form:
\[
0\leq p_1+p_2+p_3-p_{12}-p_{13}-p_{23}+p_{123}\leq1
\]
and corresponds to one of the original Bell inequalities.

If we do not want to involve all probabilities and correlations, we can use Bell-type inequalities generated by the sum of only some elementary Bell valuations on $2^N$ (for $n=3$ we can choose among $127$ inequalities). If only one of these inequalities is violated, the system is not classical.

Eventually, we will classify the Bell-like inequalities mentioned in Section~4.

These inequalities are of the form
\begin{equation}\label{equ3}
0\leq p_I+p_J-p_{I\cup J}\leq1
\end{equation}
with $I,J\subseteq N$, $I\not\subseteq J$ and $J\not\subseteq I$. Such an inequality corresponds to $f\in\mathbb R^{2^N\setminus\{\emptyset\}}$ with
\[
f(K)=\left\{
\begin{array}{rl}
 1 & \text{if }K\in\{I,J\} \\
-1 & \text{if }K=I\cup J \\
 0 & \text{otherwise}
\end{array}
\right.
\]
($K\in2^N\setminus\{\emptyset\}$). Hence
\[
g_f(K)=\sum_{\emptyset\neq L\subseteq K}f(L)=\left\{
\begin{array}{rl}
0 & \text{if }K\not\supseteq I\text{ and }K\not\supseteq J \\
1 & \text{otherwise.}
\end{array}
\right.
\]
This shows $f\in B_n$ according to Corollary~\ref{cor1}. Therefore inequalities of the form (\ref{equ3}), in particular all inequalities stated in Theorems~\ref{P2}, \ref{P3} and \ref{P4}, are all in line with our findings.

Authors' address:

Dietmar Dorninger \\
TU Wien \\
Faculty of Mathematics and Geoinformation \\
Institute of Discrete Mathematics and Geometry \\
Wiedner Hauptstra\ss e 8-10 \\
1040 Vienna \\
Austria

Helmut L\"anger \\
TU Wien \\
Faculty of Mathematics and Geoinformation \\
Institute of Discrete Mathematics and Geometry \\
Wiedner Hauptstra\ss e 8-10 \\
1040 Vienna \\
Austria, and \\
Palack\'y University Olomouc \\
Faculty of Science \\
Department of Algebra and Geometry \\
17.\ listopadu 12 \\
771 46 Olomouc \\
Czech Republic \\
helmut.laenger@tuwien.ac.at

Maciej J.\ M\c aczy\'nski \\
Uniwersytet w Bia\l ymstoku \\
Wydzia\l{} Matematyki i Informatyki \\
ul.\ Konstantego Cio\l kowskiego 1M \\
15-245 Bia\l ystok \\
Poland \\
mamacz@mini.pw.edu.pl
\end{document}